# A Conjecture about the Density of Prime Numbers.

**Lincoln Almir AMARANTE RIBEIRO**
Departamento de Física, Universidade Federal de Minas Gerais, Brasil
almir.bh@terra.com.br

**Abstract**. We present in this work a heuristic expression for the density of prime numbers. Our expression leads to results which possesses approximately the same precision of the Riemann's function in the domain that goes from 2 to $10^{10}$ at least. Instead of using a constant as was done by Legendre and others in the formula of Gauss, we try to adjust the data through a function. This function has the remarkable property: its points of discontinuity are the prime numbers.

The Prime Numbers, are natural numbers, like 2, 3, 5, 7, 11, . . . , which are not multiples of any smaller natural number, except 1. The aim of this paper is to present a heuristically derived correction to the Gauss formula to the density of Prime Numbers. Given a real number x, we can define the Density of prime numbers $\Pi(x)$, (a notation apparently due to Gauss), as the function that counts the number of prime numbers smaller then x. Before we present our conjecture it is appropiate to outline previous works on the matter.

Prime numbers were first studied by Greek mathematicians in antiquity. The Greeks had a great curiosity about the odd properties of these numbers. Soon then discovered the Fundamental Theorem of Arithmetic, which states that every integer
n >1 can be represented as a product of primes in only one way, apart from the order of the factors.
Euclid[1] (300 A.C.) in his "Elements", Book IX, proposition 20, stated that "There are infinitely many prime numbers."

Eratosthenes (275-194 B.C.) devised a procedure to discover prime numbers. It is a possible search algorithm for primes: "the sieve of Eratosthenes", which seeks to find all primes n ≤ x. It starts from all primes smaller than √x and eliminates from all integers smaller or equal to x the integer multiples of these primes. The remaining integers are the primes.

In 1737 Euler[2] demonstrated that:

$$\sum_{n=1}^{\infty} \frac{1}{n^s} = \prod_{p} \frac{1}{1-p^{-s}},$$

where the product on the right goes over all prime numbers. Euler used his formula to give an alternative proof that there are infinitely many primes. If we take $s = 1$, we get

$$\sum_{n=1}^{\infty} \frac{1}{n} = \prod_{p} \frac{1}{1-\frac{1}{p}},$$

the sum on the left is the harmonic series, which diverges. If the number of primes were finite, then the product would evaluate to a finite number. This contradiction shows that there are infinitely many primes, a new proof of Euclid´s Theorem.

Gauss[3] in his childhood studied the function $\Pi(x)$, the density of Prime Numbers. At only 15 years old, he obtained the relation between a given number and the number of prime numbers smaller than it. Through the observation that the quotient between x and log x augmented 2.3 approximately between a potency of 10 and the next potency, he conjectured:

$$\Pi(x) \approx \frac{x}{\log x},$$

that is, the limit of the quotient

$$\frac{\Pi(x)}{\dfrac{x}{\log x}},$$

as $x \to \infty$, exists and equals 1. This statement is known as the Prime Number Theorem.
In a paper from 1808, Legendre[4] claimed that $\Pi(x)$ was approximately equal to:

$$\Pi(x) = \frac{x}{\log x - B},$$

for x sufficiently large. Where B is a numerical constant whose value was suggested by him to be 1.80366.

In 1871, Gauss[5] noted that the logarithmic sum:

$$L_S(x) = 1/\log 2 + 1/\log 3 + 1/\log 4 + \ldots,$$

was a good approximation to $\Pi(x)$. Generalizing this concept he introduced the function Logarithmic Integral:

$$\text{Li}(x) = \int_2^\infty dt/\log t.$$

Gauss then proposed that:
$$\Pi(x) \approx \text{Li}(x).$$

This is a better approximation to the Density of Primes function. Gauss's assertion was completely proved in 1896 by J. Hadamard[6] and independently C. de la Vallée Poussin[7], using Riemann's work relating $\Pi(x)$ to the complex zeta function. de la Vallée Poussin also proved that Gauss' Li(x) is a better approximation to $\Pi(x)$ than x/(log x -b) no matter

what value is assigned to the constant b (and also that the best value for b is 1). They actually showed

$$\Pi(x) = \text{Li}(x) + O(xe^{-b\sqrt{\log x}}),$$

for some positive constant b. The error term depended on what was known about the zero-free region of the Riemann zeta function within the critical strip. As our knowledge of the size of this region increases, the error term decreases. Previously, P. L. Chebyshev[8] had shown that: a) There exist explicitly computable positive constants a and b such that for all x,

$$ax/\log x < \Pi(x) < bx/\log x.$$

And b) If the limit of

$$\Pi(x)/(\log x/x),$$

as $x \to \infty$ exists, it necessarily equals 1.

A study by Riemann[8] leads to the conclusion that if we consider not only the primes smaller than x, but also the primes smaller than the $n^{th}$ root of x we can obtain a better approximation to Li(x):

$$\text{Li}(x) \approx \Pi(x) + (1/2)\,\Pi(\sqrt{x}) + (1/3)\,\Pi(\sqrt[3]{x}) + ....,$$

then

$$\Pi(x) \approx \text{Li}(x) - (1/2)\,\Pi(\sqrt{x}) - (1/3)\,\Pi(\sqrt[3]{x}) + .....$$

The function of the second member is known by Riemann Function R(x). It should be pointed out that the explicit definition for the the function R(x) is

$$R(x) = \sum_{n=1}^{\infty} \frac{\mu(n)}{n} Li(x^{1/n}),$$

where μ(n) are the Möbius numbers. These are defined to be zero when *n* is divisible by a square, and otherwise to equal $(-1)^k$ where *k* is the number of distinct prime factors in *n*. As 1 has no prime factors, it follows that μ(1) = 1.

We propose the conjecture that the density of prime numbers may be written:

$$\Pi(x) = \frac{x}{\log x - f(x)},$$

and propose to search the form of this function. Using the known values of $\Pi(x)$, we obtain

$$f(x) = \log x - \frac{x}{\Pi(x)}.$$

An important property of this function is that its points of discontinuity are the prime numbers as can be seen in the figure below:

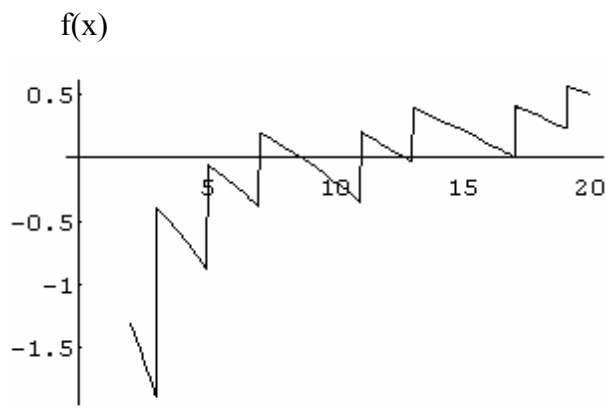

Figure 1 – Plot of x versus f(x), for x from 2 to 20. The points of discontinuity are the prime numbers

This property is obviously expected from the definition of $\Pi(x)$. If the form of f(x) function could be obtained exactly we would have a formula to calculate all the prime numbers. Unfortunately until now, it has not been possible to obtain such an exact expression. Thus we will try to obtain an approximate expression to this function, by a heuristic method. To this aim we calculate the value of f(x) for x from 10 to $10^{22}$. This last number is the last exact potency of ten for which we have a trustful exact value for $\Pi(x)$. The data are shown in table I.

| x | f(x) |
|---|---|
| 10 | 0.19741491 |
| $10^2$ | 0.60517019 |
| $10^3$ | 0.95537433 |
| $10^4$ | 1.07364387 |
| $10^5$ | 1.08757100 |
| $10^6$ | 1.07633249 |
| $10^7$ | 1.07097559 |
| $10^8$ | 1.06395401 |
| $10^9$ | 1.05662871 |
| $10^{10}$ | 1.05036512 |
| $10^{11}$ | 1.04512641 |
| $10^{12}$ | 1.04087169 |
| $10^{13}$ | 1.03734543 |
| $10^{14}$ | 1.03437617 |
| $10^{15}$ | 1.03184411 |
| $10^{16}$ | 1.02966040 |
| $10^{17}$ | 1.02775775 |
| $10^{18}$ | 1.02608510 |
| $10^{19}$ | 1.02460311 |
| $10^{20}$ | 1.02328086 |
| $10^{21}$ | 1.02209379 |
| $10^{22}$ | 1.02102214 |

Table I. values of x versus f(x)

Executing a non-linear fit to the data of table I, we obtain:

$$f(x) \approx 0.7013/x - 4.964 \exp(-0.9677x) + 0.98.$$

To simplify the calculations we used instead the variable x, the variable $\log_{10}x$. With this remark we can write:

$$\Pi(x) \approx \frac{x}{\log x \; 0.7013/\log 10x + 4.964\exp(-0.9677.\log 10x) - 0.98}$$

In table II and III below we present the comparison of the results of the values of $\Pi(x)$, for the various formulae.

For numbers smaller then 1000:

| x | Exact | This work | R(x) | Li(x) | $\frac{x}{\log x}$ |
|---|---|---|---|---|---|
| 5 | 2 | 2 | 3 | 4 | 2 |
| 10 | 4 | 4 | 4 | 6 | 4 |
| 20 | 8 | 7 | 7 | 10 | 7 |
| 30 | 10 | 10 | 10 | 13 | 9 |
| 40 | 12 | 12 | 13 | 16 | 11 |
| 50 | 15 | 14 | 15 | 18 | 13 |
| 60 | 17 | 17 | 17 | 21 | 15 |
| 70 | 19 | 19 | 19 | 23 | 16 |
| 80 | 22 | 21 | 21 | 26 | 18 |
| 90 | 24 | 23 | 24 | 28 | 20 |
| 100 | 25 | 25 | 26 | 30 | 22 |
| 200 | 46 | 44 | 45 | 50 | 38 |
| 300 | 62 | 61 | 62 | 59 | 53 |
| 400 | 78 | 78 | 78 | 85 | 67 |
| 500 | 101 | 101 | 102 | 101 | 80 |
| 600 | 109 | 109 | 110 | 118 | 94 |
| 700 | 125 | 124 | 125 | 133 | 107 |
| 800 | 139 | 139 | 140 | 148 | 120 |
| 900 | 154 | 153 | 154 | 163 | 132 |
| 1000 | 168 | 168 | 168 | 178 | 145 |

Table II – Values of the various approximations, x from 5 t0 1000

And for numbers from 10 to $10^{10}$,

| X | Exact | This work | R(x) | Li(x) | x/logx | x/(logx-1.80366) |
|---|---|---|---|---|---|---|
| 10 | 4 | 4 | 4 | 6 | 4 | 20 |
| $10^2$ | 25 | 25 | 26 | 30 | 22 | 37 |
| $10^3$ | 168 | 168 | 168 | 178 | 145 | 196 |
| $10^4$ | 1229 | 1226 | 1227 | 1246 | 1086 | 1350 |
| $10^5$ | 9592 | 9586 | 9587 | 9630 | 8686 | 10299 |
| $10^6$ | 78498 | 78533 | 78527 | 78628 | 72382 | 83251 |
| $10^7$ | 664579 | 664735 | 664667 | 664918 | 620421 | 698595 |
| $10^8$ | 5761455 | 5760802 | 5761552 | 5762209 | 5428681 | 6017926 |
| $10^9$ | 50847534 | 50848760 | 50847455 | 50849235 | 48254942 | 52855223 |
| $10^{10}$ | 455052511 | 455041196 | 455050683 | 455055614 | 434294481 | 471204883 |

Table II – Values of the various approximations, x from 10 to $10^{10}$

**Conclusion**

We made the conjecture that the density of primes is of the form:

$$\Pi(x) = \frac{x}{\log x - f(x)},$$

where f(x) is a real function of the real variable x. This function behaves at infinity as:

$$\lim_{x \to \infty} f(x) = const,$$

and so then the Prime number theorem is satisfied. The function has the amazing property: its discontinuity points occur in the prime numbers. Unfortunately we could not obtain a exact explicit expression for this function as expected. Thus we are compelled to appeal to a heuristic approach. In this manner we tried to fit a non-linear function to the data we obtained by considering each point as the quantity needed to make the expression:

$$\Pi(x) \approx \frac{x}{\log x - f(x)},$$

as precise as possible. As a result we obtain an expression for the function that when inserted in the formula of Gauss gives results comparable to the Riemann function in the domain 2 to $10^{10}$ at least. One advantage is that we can obtain $\Pi(x)$ using a scientific pocket calculator such, as an HP-15C, with a great speed and ease. We hope that mathematicians moved by this form of the function encounter clues that will enable the exact functional form of f(x) to be determined if possible.

**References**

[1]- Euclid, " Euclid's Elements: all thirteen books complete in one volume". The Thomas L. Heath Translation. Dana Densmore, Editor. Santa Fe, New Mexico: Green Lion Press, 2002, 529 pp.

[2] - Euler, Leonhard, Variae observations circa series infinitas, *Commentarii academiae scientiarumPetropolitanae* **9** (1737), 1744, p. 160-188. Reprinted in *Opera Omnia* Series I volume 14, pps. 216-244.

[3] –This result was unpublished. Gauss wrote a letter to Johann Lencke at December, 25, 1849. In this letter he reports the connection between prime numbers and logarithms and presents his formula.

[4] - Legendre, A. M. "Essai sur la theorie des nombres", Paris: Duprat, 1798

[5] - Gauss, C. F. *Werke, Band 10, Teil I.* 1863, p. 10.

[6] - Hadamard, J, 1896, "Sur la distribuition des zeros de la function ζ(s) et ces consequences aritmetiques. Bull. Soc. Math. France, XXIV, 199-220

[7] – de la Vallee Poussin C.J., Recherches analytiques sur la theorie des nombres premiers, Ann. Soc. Sci. Bruxelles 20 (1896), pps. 183-256.

[8] - Chebyshev, P. L. "Mémoir sur les nombres premiers." *J. math. pures appl.* **17**, 1852.

[9] - B. Riemann, ¨Uber die Anzahl der Primzahlen unter einer gegebenen Größe, Monatsberichte der Berliner Akademie, 1859, pps. 671-680